\documentstyle[11pt,leqno]{article}
\newtheorem{theorem}{{\sc Theorem}}
\newcommand{\bt}{\begin{theorem}}
\newcommand{\et}{\end{theorem}}
\newcommand{\newsection}[1]{\setcounter{equation}{0} \setcounter{theorem}{0}
\section{#1}}

\newcommand{\NI}{\noindent}
\newcommand{\bea}{\begin{eqnarray}}
\newcommand{\eea}{\end{eqnarray}}

\def \spec#1 {\mathop{#1}}

\def \b #1 {\bf #1}
\newcommand {\nnb }{\nonumber}

\newcommand {\CC}{\centerline}

\newcommand{\clf}{{\cal F}}

\newcommand{\ity}{\infty}
\newcommand{\raro}{\rightarrow}

\newcommand{\vsp}{\vskip 1em}
\newcommand{\vspp}{\vskip 2em}

\newcommand{\ve}{\varepsilon}

\newcommand{\be}{\begin{equation}}
\newcommand{\ee}{\end{equation}}
\newcommand{\ben}{\begin{eqnarray*}}
\newcommand{\een}{\end{eqnarray*}}

\begin{document}

\sloppy
\CC {\Large{\bf Berry-Esseen Type Bound for Fractional Ornstein-Uhlenbeck}}
 \CC {\Large {\bf Type Process Driven by Sub-Fractional Brownian Motion }}
\vsp
\CC {\bf B.L.S. PRAKASA RAO}
\CC {\bf CR Rao Advanced Institute of Research in Mathematics,}
\CC{\bf  Statistics and Computer Science, Hyderabad 500046, India}
\vsp
\NI{\bf Abstract:} We obtain a Berry-Esseen type bound for the distribution of the maximum likelihood estimator
of the drift parameter  for fractional Ornstein-Uhlenbeck type process driven by sub-fractional Brownian motion.
\vsp
\NI{\bf Keywords and phrases}: Fractional Ornstein-Uhlenbeck type process ; sub-fractional Brownian motion; Maximum likelihood estimation; Berry-Esseen type bound.
\vspp
\NI {\bf AMS Subject classification (2010): } Primary 62M09, Secondary 60G22.
\newsection{Introduction}

Statistical inference  for fractional diffusion processes satisfying stochastic
differential equations driven by a fractional Brownian motion (fBm)  has been studied
earlier and a comprehensive survey of various methods is given in Prakasa
Rao (2010). There has been a recent interest to study similar problems for
stochastic processes driven by a sub-fractional Brownian motion. Bojdecki et al. (2004) introduced a centered Gaussian process $\zeta^H= \{\zeta^H(t), t\geq 0\}$ called  {\it sub-fractional Brownian motion} (sub-fBm) with the covariance function
$$C_H(s,t)= s^{2H}+t^{2H}-\frac{1}{2}[(s+t)^{2H}+|s-t|^{2H}]$$
where $0<H<1.$ The increments of this process are not stationary and are more weakly correlated on non-overlapping intervals than those of a fBm. Tudor (2009) introduced a Wiener integral with respect to a sub-fBm. Tudor ( 2007 a,b, 2008, 2009) discussed some properties related to sub-fBm and its corresponding stochastic calculus. By using a fundamental martingale associated to sub-fBm, a Girsanov type theorem is obtained in Tudor (2009). Diedhiou et al. (2011) investigated parametric estimation for a stochastic differential equation (SDE) driven by a sub-fBm. Mendy (2013) studied parameter estimation for the sub-fractional Ornstein-Uhlenbeck process defined by the stochastic differential equation
$$dX_t=\theta X_tdt+d\zeta^H(t), t \geq 0$$ where $H>\frac{1}{2}.$
This is an analogue of the
Ornstein-Uhlenbeck process, that is, a continuous time first order
autoregressive process $X=\{X_t, t \geq 0\}$ which is the solution of a
one-dimensional homogeneous linear stochastic differential equation driven
by a sub-fBm $\zeta^H= \{\zeta_t^H, t \geq 0\}$ with
Hurst parameter $H.$ Mendy (2013) proved that the least squares estimator estimator $\tilde \theta_T$ is strongly consistent as $T \raro \ity.$
Kuang and Xie (2015) studied properties of maximum likelihood estimator for sub-fBm through approximation by a random walk. Kuang and Liu (2015) discussed about the $L^2$-consistency  and strong consistency of the maximum likelihood estimators for the sub-fBm with drift based on discrete observations. Yan et al. (2011) obtained the Ito's formula for sub-fractional Brownian motion with Hurst index $H>\frac{1}{2}.$  Shen and Yan (2014) studied estimation for the drift of sub-fractional Brownian motion and constructed a class of biased estimators of James-Stein type which dominate the maximum likelihood estimator under the quadratic risk. El Machkouri et al. (2016) investigated the asymptotic properties of the least squares estimator for non-ergodic Ornstein-Uhlenbeck process driven by Gaussian processes, in particular, sub-fractional Brownian motion. In a recent paper, we have investigated optimal estimation of a signal perturbed by a sub-fractional Brownian motion in Prakasa Rao (2017b). Some maximal and integral inequalities for a sub-fBm were derived in Prakasa Rao (2017a).  Parametric estimation for linear stochastic differential equations driven by a sub-fractional Brownian motion is studied in Prakasa Rao (2017c). We now obtain a Berry-Esseen type bound for the distribution of the maximum likelihood estimator for the drift parameter of a fractional Ornstein-Uhlenbeck type process driven by a sub-fractional Brownian motion.

\newsection{Preliminaries}
Let $(\Omega, \clf, (\clf_t), P) $ be a stochastic basis satisfying the
usual conditions and the processes discussed in the following are
$(\clf_t)$-adapted. Further the natural filtration of a process is
understood as the $P$-completion of the filtration generated by this
process.

Let $\zeta^H= \{\zeta_t^H, t \geq 0 \}$ be a  normalized  {\it sub-fractional Brownian
motion} (sub-fBm)  with Hurst parameter $H \in (0,1)$, that is, a Gaussian process
with continuous sample paths such that $\zeta_0^H=0, E(\zeta_t^H)=0$ and
\be
E(\zeta_s^H \zeta_t^H)= t^{2H}+s^{2H}-\frac{1}{2}[(s+t)^{2H}+|s-t|^{2H}], t \geq 0, s \geq 0.
\ee

Bojdecki et al. (2004) noted that the process
$$\frac{1}{\sqrt{2}}[W^H(t)+W^H(-t)], t \geq 0,$$
where $\{W^H(t), -\ity<t<\ity\}$ is a fBm, is a centered Gaussian process with the same covariance function as that of a  sub-fBm. This proves the existence of a sub-fBm. They proved the following result concerning properties of a sub-fBm.
\vsp
\NI{\bf Theorem 2.1:} {\it Let $\zeta^H= \{\zeta^H(t), t \geq 0\}$ be a sub-fBm defined on a filtered probability space $(\Omega, \clf, (\clf_t,t\geq 0),P).$ Then the following properties hold.

(i) The process $\zeta^H$ is self-similar, that is, for every $a>0,$
$$\{\zeta^H(at),t\geq 0\} \stackrel {\Delta}=  \{a^H\zeta^H(t), t\geq 0\}$$
in the sense that the processes, on both the sides of the equality sign, have the same finite dimensional distributions.

(ii) The process $\zeta^H$ is not Markov and it is not a semi-martingale.

(iii) Let $R_H(s,t)$ denote the covariance function of standard fractional Brownian motion with Hurst index $H.$ For all $s,t \geq 0,$ the covariance function $C_H(s,t)$ of the process $\zeta^H$ is positive for all $s>0,t>0.$ Furthermore

$$C_H(s,t)>R_H(s,t) \;\;\mbox{if}\;\;H<\frac{1}{2}$$
and
$$C_H(s,t)< R_H(s,t) \;\;\mbox{if}\;\;H>\frac{1}{2}.$$

(iv) Let $\beta_H= 2-2^{2H-1}.$ For all $s\geq 0, t\geq 0,$

$$\beta_H(t-s)^{2H}\leq E[\zeta^H(t)-\zeta^H(s)]^2\leq (t-s)^{2H},\;\;\mbox{if}\;\;H>\frac{1}{2}$$
and
$$(t-s)^{2H}\leq E[\zeta^H(t)-\zeta^H(s)]^2\leq \beta_H (t-s)^{2H},\;\;\mbox{if}\;\;H<\frac{1}{2}$$
and the constants in the above inequalities are sharp.

(v) The process $\zeta^H $ has continuous sample paths almost surely and, for each $0<\epsilon<H$ and $T>0,$  there exists a random variable $K_{\epsilon,T}$ such that}
$$|\zeta^H(t)-\zeta^H(s)|\leq K_{\epsilon,T} |t-s|^{H-\epsilon}, 0\leq s,t \leq T.$$
\vsp
Let $f:[0,T]\raro R$ be a measurable  function and $\alpha  >0,$ and  $\sigma$ and $\eta$ be real. Define the Erdeyli-Kober-type fractional integral
\be
(I_{T,\sigma,\eta}^\alpha f)(s)=\frac{\sigma s^{\alpha \eta}}{\Gamma(\alpha)}\int_s^T\frac{t^{\sigma(1-\alpha-\eta)-1}f(t)}{(t^\sigma-s^\sigma)^{1-\alpha}}dt, s\in [0,T],
\ee
and the function
\bea
n_H(t,s)& =& \frac{\sqrt{\pi}}{2^{H-\frac{1}{2}}}I_{T,2,\frac{3-2H}{4}}^{H-\frac{1}{2}}(u^{H-\frac{1}{2}})I_{[0,t)}(s)\\\nonumber
&=& \frac{2^{1-H}\sqrt{\pi}}{\Gamma(H-\frac{1}{2})}s^{\frac{3}{2}-H}\int_0^t(x^2-s^2)^{H-\frac{3}{2}}dx \;I_{(0,t)}(s).\\\nonumber
\eea
The following theorem is due to Dzhaparidze and Van Zanten (2004) (cf. Tudor (2009)).
\vsp
\NI{\bf Theorem 2.2:} {\it The following representation holds, in distribution, for a sub-fBm $\zeta^H$:
\be
\zeta^H_t \stackrel {\Delta} = c_H\int_0^t n_H(t,s)dW_s, 0\leq t \leq T
\ee
where
\be
c_H^2= \frac{\Gamma (2H+1)\;\sin(\pi H)}{\pi}
\ee
and $\{W_t,t\geq 0\}$ is the standard Brownian motion.}
\vsp
Tudor (2009) has defined integration of a non-random function $f(t)$ with respect to a sub-fBm $\zeta^H$ on an interval $[0,T]$ and obtained a representation of this integral as a Wiener integral for a suitable transformed function $\phi_f(t)$ depending on $H$ and $T.$ For details, see Theorem 3.2 in Tudor (2009).
\vsp
Tudor (2007b) obtained the prediction formula for a sub-fBm. For any $0<H<1,$ and $0<a<t,$
\be
E[\zeta_t^H|\zeta_s^H, 0 \leq s \leq a]=\zeta_a^H+ \int_0^a \psi_{a,t}(u)d \zeta_u^H
\ee
where
\be
\psi_{a,t}(u)= \frac{2\; \sin(\pi (H-\frac{1}{2}))}{\pi} u(a^2-u^2)^{\frac{1}{2}-H}\int_a^t\frac{(z^2-a^2)^{H-\frac{1}{2}}}{z^2-u^2}z^{H-\frac{1}{2}}dz.
\ee
Let
\be
M_t^H=d_H \int_0^ts^{\frac{1}{2}-H}dW_s = \int_0^tk_H(t,s)d\zeta_s^H
\ee
where
\be
d_H=\frac{2^{H-\frac{1}{2}}}{c_H \Gamma(\frac{3}{2}-H)\sqrt{\pi}},
\ee
\be
k_H(t,s)=d_Hs^{\frac{1}{2}-H}\psi_H(t,s),
\ee
and
\ben
\psi_H(t,s)&=& \frac{s^{H-\frac{1}{2}}}{\Gamma(\frac{3}{2}-H)}[t^{H-\frac{3}{2}}(t^2-s^2)^{\frac{1}{2}-H}-\\\nonumber
&&\;\;\;\; (H-\frac{3}{2})\int_s^t(x^2-s^2)^{\frac{1}{2}-H}x^{H-\frac{3}{2}}dx]I_{(0,t)}(s).\\\nonumber
\een
\vsp
It can be shown that the process $ M^H=\{M^H_t, t \geq 0\}$ is a Gaussian martingale (cf. Tudor (2009), Diedhiou et al. (2011)) and is called the {\it sub-fractional fundamental martingale}. The filtration generated by this martingale is the same as the filtration $\{\clf_t, t\geq 0\}$ generated by the sub-fBm $\zeta^H$ and the quadratic variation $<M^H>_s$ of the martingale $M^H$ over the interval $[0,s]$ is equal to $w_s^H=\frac{d_H^2}{2-2H}s^{2-2H}= \lambda_H s^{2-2H}$ (say). For any measurable function $f:[0,T] \rightarrow R$ with $\int_0^Tf^2(s)s^{1-2H} ds<\ity,$ define the probability measure $Q_f$ by
\ben
\frac{dQ_f}{dP}|_{\clf_t}&=& \exp(\int_0^tf(s)dM_s^H-\frac{1}{2}\int_0^tf^2(s)d<M^H>(s))\\
&=&\exp(\int_0^tf(s)dM_s^H-\frac{d_H^2}{2}\int_0^tf^2(s)s^{1-2H}ds)
\een
where $P$ is the underlying probability measure. Let
\be
(\psi_H f)(s)=\frac{1}{\Gamma(\frac{3}{2}-H)}I_{0,2,\frac{1}{2}-H}^{H-\frac{1}{2}}f(s)
\ee
where, for $\alpha >0,$
\be
(I_{0,\sigma,\eta}^\alpha f)(s)=\frac{\sigma s^{-\sigma(\alpha+ \eta)}}{\Gamma(\alpha)}\int_0^s \frac{t^{\sigma(1+\eta)-1}f(t)}{(t^\sigma-s^\sigma)^{1-\alpha}}dt, s\in[0,T].
\ee
Then the following Girsanov type theorem holds for the sub-fBm process (Tudor (2009)).
\vsp
\NI{\bf Theorem 2.3:} {\it The process
$$\zeta_t^H- \int_0^t(\psi_Hf)(s)ds, 0\leq t \leq T$$
is a sub-fbm with respect to the probability measure $Q_f.$ In particular, choosing the function $f\equiv a \in R$, it follows that the process $\{\zeta_t^H-at, 0\leq t \leq T\}$ is a sub-fBm under the probability measure $Q_f$ with $f \equiv a\in R.$}
\vsp
Let $Y=\{Y_t, t \geq 0\}$ be a stochastic process defined on the filtered probability space $(\Omega, \clf, (\clf_t, t\geq 0), P)$ and suppose the process $Y$ satisfies the stochastic differential equation
\be
dY_t= C(t) dt+ d\zeta_t^H, t \geq 0
\ee
where the process $\{C(t), t \geq 0\},$ adapted to the filtration $\{\clf_t, t \geq 0\},$ such that the process
\be
R_H(t)= \frac{d}{dw_t^H}\int_0^tk_H(t,s)C(s)ds, t \geq 0
\ee
is well-defined and the derivative is understood in the sense of absolute continuity with respect to the measure generated by the function $w_H.$ Differentiation with respect to $w_t^H$ is understood in the sense:
$${dw_t^H}= \lambda_H (2-2H)t^{1-2H}dt$$
and
$$\frac{df(t)}{dw_t^H}= \frac{df(t)}{dt} / \frac{dw_t^H}{dt}.$$
Suppose the process $\{R_H(t), t \geq 0\},$ defined over the interval $[0,T]$ belongs to the space $L^2([0,T],dw_t^H).$ Define
\be
\Lambda_H(t)=\exp\{\int_0^tR_H(s)dM_s^H-\frac{1}{2}\int_0^t[R_H(s)]^2dw_s^H\}
\ee
with $E[\Lambda_H(T)]=1$ and the distribution of the process $\{Y_t, 0\leq t \leq T\}$ with respect to the measure $P^Y= \Lambda_H(T)\;P$ coincides with the distribution of the process $\{\zeta_t^H, 0 \leq t \leq T\|$ with respect to the measure $P.$
\vsp
We call the process $\Lambda^H$ as the {\it likelihood process} or the Radon-Nikodym derivative $\frac{dP^Y}{dP}$ of the measure $P^Y$ with respect to the measure $P.$
\vsp
\vsp
Tudor (2009) derived the following Girsanov type formula.
\vsp
\NI{\bf Theorem 2.4:} {\it Suppose the assumptions of Theorem 2.2 hold. Define
\be
\Lambda_H(T)= \exp \{ \int_0^T R_H(t)dM_t^H - \frac{1}{2}\int_0^TR_H^2(t)d w_t^H\}.
\ee
Suppose that $E(\Lambda_H(T))=1.$ Then the measure $P^*= \Lambda_H(T) P$ is
a probability measure and the probability measure of the process $Y$ under
$P^*$ is the same as that of the process $V$ defined by}
\be
V_t= \int_0^t d\zeta_s^H, 0 \leq t \leq T.
\ee.
\newsection{Main Results}
Let us consider the stochastic differential equation
\be
dX(t)= \theta \; X(t)dt + d\zeta_t^H, X(0)=0, t \geq 0
\ee
where $\theta \in \Theta \subset R, \zeta^H=\{\zeta_t^H, t \geq 0\}$ is a
sub-fractional Brownian motion with known Hurst parameter $H.$ In other words $X=\{X(t),
t \geq 0\}$ is a stochastic process satisfying the stochastic integral
equation
\be
X(t)= \theta \int_0^tX(s)ds + \int_0^t d\zeta_s^H, t \geq 0.
\ee
We call such a process as fractional Ornstein-Uhelenbeck type process  driven by sub-fractional Brownian motion. Diedhiou et al. (2011) and Mendy (2013) investigated parametric estimation for such a stochastic differential equation driven by a sub-fBm. We will now obtain a Berry-Esseen type bound for the distribution of the maximum likelihood estimator for the drift parameter for such processes.
\vsp
Let
\be
C(\theta,t)= \theta \; X(t), t \geq 0
\ee
and assume that the sample paths of the process $\{C(\theta,t), t \geq 0\}$ are smooth enough so that
the process
\be
R_{H,\theta}(t) =\theta \frac{d}{dw_t^H}\int_0^tk_H(t,s)X(s)ds, t \geq 0
\ee
is well-defined where $w_t^H$ and $k_H(t,s)$ are as defined in Section 2. Suppose the sample paths of the process $\{R_{H,\theta}(t), 0
\leq t \leq T \}$ belong almost surely to $L^2([0,T], dw_t^H).$ Define
\be
Z_t= \int_0^t k_H(t,s) d X_s, t \geq 0.
\ee
Then the process $Z= \{Z_t,t \geq 0\}$ is an $(\clf_t)$-semimartingale
with the decomposition
\be
Z_t= \int_0^t R_{H, \theta}(s)dw_s^H + M_t^H, t \geq 0
\ee
where $M^H$ is the fundamental martingale defined by the equation (2.8) and the process
$X$ admits the representation
\be
X_t = \int_0^t K_H(t,s) dZ_s
\ee
where the function
$$
K_H(t,s)= \frac{c_H}{d_H}s^{H-\frac{1}{2}}n_H(t,s).
$$
\vsp
Let $P_\theta^T$ be the measure induced by the process $\{X_t, 0 \leq t \leq T\}$ when $\theta$
is the true parameter. Following Theorem 2.4, we get that the
Radon-Nikodym derivative of $P_\theta^T$ with respect to $P_0^T$ is given by
\be
\frac{dP_\theta^T}{dP_0^T}= \exp[ \int_0^T R_{H,
\theta}(s)dZ_s - \frac{1}{2} \int_0^T R_{H, \theta}^2(s)dw_s^H].
\ee
\vsp
\NI {\bf Maximum likelihood estimation}

We now consider the problem of estimation of the parameter $\theta$ based
on the observation of the process $ X= \{X_t , 0 \leq t \leq T\}$ and
study its asymptotic properties as $T \raro \ity.$
\vsp
\NI {\bf Strong consistency:}

Let $L_T(\theta)$ denote the Radon-Nikodym derivative
$\frac{dP_\theta^T}{dP_0^T}.$ The maximum likelihood estimator (MLE) is
defined by the relation
\be
L_T(\hat \theta_T) = \sup_{\theta \in \Theta} L_T(\theta).
\ee
We assume that there exists a measurable maximum likelihood estimator.
Sufficient conditions can be given for the existence of such an estimator
(cf. Lemma 3.1.2, Prakasa Rao (1987)).
Note that
\bea
R_{H,\theta}(t) &= &\theta \frac{d}{dw_t^H}\int_0^tk_H(t,s)X(s)ds\\\nnb
& = & \theta J(t). (say)
\eea
Then
\be
\log L_T(\theta)= \theta \int_0^T J(t) dZ_t
-\frac{1}{2}\theta^2\int_0^T J^2(t) dw_t^H
\ee
and the likelihood equation is given by
\be
\int_0^T J(t) dZ_t - \theta \int_0^T J^2(t) dw_t^H =0.
\ee
Hence the MLE $\hat \theta_T$ of $\theta$ is given by
\be
\hat \theta_T = \frac {\int_0^T J(t) dZ(t)}{\int_0^T J^2(t) dw_t^H}.
\ee
Let $\theta_0$ be the true parameter. Using the fact that
\be
dZ_t= \theta_0 J(t)) dw_t^H + dM_t^H,
\ee
it can be shown that
\be
\frac{dP_{\theta}^T}{dP_{\theta_0}^T}= \exp[ (\theta- \theta_0)\int_0^TJ(t)dM_t^H - \frac{1}{2} (\theta-\theta_0)^2 \int_0^T J^2(t) d w_t^H]
.
\ee
Following this representation of the Radon-Nikodym derivative, we obtain that
\be
\hat \theta_T- \theta_0= \frac{\int_0^T J(t) d M_t^H}{\int_0^TJ^2(t)dw_t^H}.
\ee
\vsp
We now discuss the problem of estimation of the parameter $\theta$ on the
basis of the observation of the process $X$ or equivalently the process
$Z$ on the interval $[0,T].$
\vsp
\NI {\bf Theorem 3.1:} {\it The maximum likelihood estimator $\hat \theta_T$ is
strongly consistent, that is,
\be
\hat \theta_T \raro \theta_0 \;\;\mbox{a.s}\;\; [P_{\theta_0}] \;\mbox
{as}\;\; T \raro \ity
\ee
provided}
\be
\int_0^TJ^2(t) dw_t^H \raro \ity\;\; \mbox{a.s}\;\; [P_{\theta_0}]\;\;
\mbox{as}\;\; T \raro \ity.
\ee
\vsp
\NI{\bf Proof:} This theorem follows by observing that the process
\be
\gamma_T \equiv \int_0^T J(t)dM_t^H, t \geq 0
\ee
is a local martingale with the quadratic variation process
\be
<\gamma>_T= \int_0^T J^2(t)dw_t^H
\ee
and applying the Strong law of large numbers (cf. Liptser (1980); Liptser and Shiryayev (1989); Prakasa
Rao (1999), p. 61) under the condition (3.18) stated above.
\vsp
\noindent {\bf Remark:} For the case sub-fractional Ornstein-Uhlenbeck process investigated in
Mendy (2013), it can be checked that the condition
stated in equation (3.18) holds and hence the maximum likelihood estimator
$\hat \theta_T $ is strongly consistent as $T \raro \ity.$
\vsp
\NI {\bf Limiting distribution:}

We now discuss the limiting distribution of the MLE $\hat \theta_T$ as $T \raro \ity.$
\vsp
\NI {\bf Theorem 3.2:} {\it Suppose that the process $\{\gamma_t, t \geq 0\}$ is a local continuous
martingale and that there exists a norming function $I_t, t \geq 0$ such
that
\be
I_T^2 <\gamma_T> = I_T^2 \int_0^T J^2(t)dw_t^H \raro \eta^2 \;\; \mbox{in
probability}\;\; \mbox{as}\;\; T\raro \ity
\ee
where $ I_T \raro 0 $ as $T \raro \ity $ and $\eta$ is a random variable
such that $P(\eta >0)=1.$ Then
\be
(I_T \gamma_T, I_T^2 <\gamma_T>) \raro (\eta Z, \eta^2) \mbox{ in law}\;\;
\mbox{as}\;\; T \raro \ity
\ee
where the random variable $Z$ has the standard normal distribution and the
random variables $Z$ and $\eta$ are independent.}
\vsp
\NI {\bf Proof:} This theorem follows as a consequence of the central limit
theorem for martingales (cf. Theorem 1.49 ; Remark 1.47 , Prakasa
Rao (1999), p. 65).
\vsp
Observe that
\be
I_T^{-1}(\hat \theta_T - \theta_0)  =  \frac{I_T \gamma_T}{I_T^2<\gamma_T>}
\ee
Applying the Theorem 3.2, we obtain the following result.
\vsp
\NI {\bf Theorem 3.3:} {\it Suppose the conditions stated in the Theorem 3.2
hold. Then
\be
I_T^{-1}(\hat \theta_T - \theta_0) \raro \frac{ Z}{\eta} \mbox{ in
law}\;\; \mbox{as}\;\; t \raro \ity
\ee
where the random variable $Z$ has the standard normal distribution and the
random variables $Z$ and $\eta$ are independent.}
\vsp
\NI {\bf Remarks:} If the random variable $\eta$ is a constant with probability
one, then the limiting distribution of the maximum likelihood estimator is
normal with mean 0 and variance $\eta^{-2}.$ Otherwise it is a mixture of
the normal distributions with mean zero and variance $\eta^{-2}$ with the
mixing distribution as that of $\eta.$
\vsp
\newsection{Berry-Esseen Type Bound}
Let $\theta_0$ be the true parameter. In addition to the conditions stated in Section 3, suppose that the random variable $\eta$ is a positive constant with probability one under $P_{\theta_0}$-measure. Theorem 3.3 implies that
\be
I_T^{-1}(\hat \theta_T-\theta_0) \raro N(0,\eta^{-2})\;\;\mbox{in law as }\;\;T \raro \ity
\ee
under $P_{\theta_0}$-measure  where $N(0, \sigma^2)$ denoted the Gaussian distribution with mean zero and variance $\sigma^2.$ We would now like to obtain the rate of convergence in this limit leading to a Berry-Esseen type bound.
\vsp
Suppose there exists non-random positive functions $\delta_T $ and $\epsilon_T$ decreasing to zero as $T\raro \ity$ such that
\be
\delta_T^{-1}\epsilon^2_T \raro \ity \;\;\mbox{as}\;\;T\raro \ity
\ee
and
\be
\sup_{\theta \in \Theta}P_{\theta}^T(|\delta_T <\gamma>_T-1|\geq \epsilon_T)=O(\epsilon_T^{1/2})
\ee
where the process $\{\gamma_T, T \geq 0\}$ is as defined by equation (3.19).  Note that the process $\{\gamma_T, T \geq 0\}$ is a locally square integrable continuous martingale. From the results on the representation of locally square integrable continuous martingales (cf. Ikeda and Watanabe (1981), Chapter II, Theorem 7.2), it follows that there exists a standard
Wiener process $\{B(t), t \geq 0\}$ adapted to $(\clf_t)$
such that $\gamma_t= B(<\gamma>_T),  t \geq 0.$ In particular \be \gamma_T
\delta_T^{1/2} = B(<\gamma>_T \delta_T)\;\;
\mbox{a.s.}\;\;[P_{\theta_0}] \ee for all $T \geq 0.$
\vsp We use the following lemmas in the sequel.
\vsp
\NI {\bf Lemma 4.1:} {\it Let $(\Omega, \clf, P)$ be a probability space and
$f$ and $g$ be $\clf$-measurable functions. Then, for any $\ve >0,$
\bea
\lefteqn{\sup_x|P(\omega: \frac{f(\omega)}{g(\omega)} \leq x)- \Phi(x)|}
\\ \nnb
& & \;\;  \leq \sup_y|P(\omega: f(\omega)  \leq y)- \Phi(x)| + P(\omega:
|g(\omega) -1| > \ve) + \ve
\eea
where $\Phi(x)$ is the distribution function of the standard  Gaussian
distribution.}

\NI {Proof:} See Michael and Pfanzagl (1971).
\vsp
\NI {\bf Lemma 4.2:} {\it Let $\{B(t), t \geq 0\}$ be a standard Wiener process and
$V$ be a nonnegative random variable. Then, for every $x \in R$
and $\ve > 0,$}
\be |P(B(V) \leq x) - \Phi(x)| \leq (2\ve)^{1/2} +
P(|V-1| > \ve).
\ee
\NI {Proof:} See Hall and Heyde (1980), p.85.
\vsp
Let us fix $\theta \in \Theta.$ It is clear from the earlier remarks
that
\be
\gamma_T= <\gamma>_T I_T^{-1}(\hat \theta_T - \theta)
\ee
under $P_{\theta}$-measure. Then it follows, from the Lemmas 4.1 and
4.2, that
\bea
\lefteqn{P_\theta[\delta_T^{-1/2}I_T^{-1}(\hat \theta_T - \theta) \leq x]- \Phi(x)|} \\ \nnb
& &  = |P_\theta[\frac {\gamma_T}{<\gamma>_T}\delta_T^{-1/2} \leq x]- \Phi(x)| \\ \nnb
& & = |P_\theta[\frac{\gamma_T/\delta_T^{-1/2}}{<\gamma>_T/ \delta_T^{-1}} \leq x]-
\Phi(x)| \\ \nnb
& & \leq  \sup_x |P_\theta[\gamma_T \delta_T^{1/2} \leq x]- \Phi(x)| \\ \nnb
& & \;\; \; \; + P_\theta[|\delta_T <\gamma>_T - 1| \geq \ve_T] + \ve_T \\ \nnb
& & =\sup_y |P(B(<\gamma>_T \delta_T) \leq y) - \Phi(y)| +
P_\theta[|\delta_T <\gamma>_T - 1| \geq \ve_T] + \ve_T \\ \nnb
& & \leq (2 \ve_T)^{1/2} + 2 P_\theta [|\delta_T <\gamma>_T - 1| \geq \ve_T] +
\ve_T.
\eea
It is clear that the bound obtained above is of the order
$O(\ve_T^{1/2})$
under the condition (4.3) and it is uniform in $\theta \in \Theta.$
Hence we have the following result giving a Berry-Esseen type bound for the distribution of the MLE.
\vsp
\NI{ \bf Theorem 4.3:} {\it Under the conditions (4.2) and (4.3),}
\bea
\lefteqn{\sup_{\theta \in \Theta} \sup_{x \in R}|P_\theta[\delta_T^{-1/2}I_T^{-1}(\hat\theta_T - \theta) \leq x]- \Phi(x)|} \\ \nnb
& & \;\; \leq (2 \ve_T)^{1/2} + 2 P_\theta[|\delta_T <\gamma>_T - 1| \geq \ve_T] +
\ve_T = O(\ve_T^{1/2}).
\eea
\vsp
As a consequence of this result, we have the following
theorem giving the rate of convergence of the  MLE $\hat \theta_T.$
\vsp
\NI {\bf Theorem 4.4:} {\it  Suppose the conditions (4.2) and (4.3) hold. Then there exists a constant $c >0$ such
that  for every $d >0,$}
\be \sup_{\theta \in \Theta}P_\theta[I_T^{-1}|\hat \theta_T - \theta| \geq d] \leq c
\ve_T^{1/2} + 2 P_\theta[|\delta_T <\gamma>_T - 1| \geq \ve_T] =
O(\ve_T^{1/2}).
\ee
\vsp
\NI {Proof:}  Observe that
\bea
\lefteqn{\sup_{\theta \in \Theta} P_\theta[I_T^{-1}|\hat\theta_T - \theta| \geq d]} \\ \nnb
& & \leq \sup_{\theta \in \Theta} |P_\theta[\delta_T^{-1/2}I_T^{-1}(\hat \theta_T -\theta) \geq d \delta_T^{-1/2}] - 2(1- \Phi(d \delta_T^{-1/2}))|\\
\nnb
& & \;\;\;\; + 2 (1- \Phi(d\delta_T^{-1/2})) \\ \nnb
& & \leq (2 \ve_T)^{1/2} + 2  \sup_{\theta \in \Theta} P_\theta[|\delta_T<\gamma>_T - 1| \geq \ve_T] +  \ve_T \\ \nnb
 & & \;\;\; \; + 2 d^{-1}\delta_T^{1/2} (2\pi)^{-1/2} \exp[-\frac{1}{2}\delta_T^{-1}d^2]
\eea
by Theorem 4.3 and the inequality
\be 1-\Phi(x) < \frac{1}{x\sqrt{2\pi}}\exp[-\frac{1}{2}x^2]
\ee
for all $x >0$ (cf. Feller (1968), p.175). Since
$$\delta_T^{-1} \ve^2_T \raro \ity \;\; \mbox{as}\;\;T \raro \ity  $$
by the condition (4.2), it follows that
\be
\sup_{\theta \in \Theta} P_\theta[I_T^{-1}|\hat \theta_T - \theta| \geq d] \leq c  \ve_T^{1/2} + 2  \sup_{\theta \in \Theta} P_\theta[|\delta_T<R>_T - 1| \geq \ve_T]
\ee
for some constant $c > 0$ and the last term is of the order
$O(\ve_T^{1/2})$ by the condition (4.3). This proves Theorem 4.4.
\vsp
\NI{\bf Acknowledgement:} This work was supported under the scheme ``INSA Senior Scientist" of the Indian National Science Academy at the CR Rao Advanced Institute of Mathematics, Statistics and Computer Science, Hyderabad 500046, India.
\vsp
\NI {\bf References}
\begin{description}

\item Artemov, A.V. and Burnaev, E.V. (2016) Optimal estimation of a signal perturbed by a fractional Brownian motion, {\it Theory Probab. Appl.}, {\bf 60}, 126-134.

\item Bojdecki, T. and Gorostiza, A. Talarczyk (2004) Sub-fractional Brownian motion and its relation to occupation times, {\it Statist. Probab. Lett.}, {\bf 69}, 405-419.

\item Diedhiou, A., Manga, C. and Mendy, I. (2011) Parametric estimation for SDEs with additive sub-fractional Brownian motion, {\it Journal of Numerical Mathematics and Stochastics}, {\bf 3}, 37-45.

\item Dzhaparidze, K. and Van Zanten, H. (2004) A series expansion of fractional Brownian motion, {\it Probab. Theory Related Fields}, {\bf 103}, 39-55.

\item El Machkouri, Mohammed; Es-Sebaiy, Khalifa and Ouknine, Youssef (2016) Least squares estimation for non-ergodic Ornstein-Uhlenbeck processes driven by Gaussian processes, arXiv:1507.00802v2 [math.PR]  26 Sep 2016.

\item Feller, W. (1968) {\it An Introduction to Probability Theory and its Applications}, New York: Wiley.

\item Hall, P. and Heyde, C.C. (1980) {\it Martingale Limit Theory and its Applications}, New York: Academic Press.

\item Ikeda, N. and Watanabe, S. (1981) {\it Stochastic Differential Equations and Diffusion Processes}, Amsterdam: North-Holland.

\item Kuang, Nenghui and Liu, Bingquan (2015) Parameter estimations for the sub-fractional Brownian motion with drift at discrete observation, {\it Brazilian Journal of Probability and Statistics}, {\bf 29}, 778-789.

\item Kuang, Nenghui and Xie, Huantin (2015) Maximum likelihood estimator for the sub-fractional Brownian motion approximated by a random walk, {\it Ann. Inst. Statist. Math.}, {\bf 67}, 75-91.

\item Liptser, R.S. (1980) A strong law of large numbers, {\it Stochastics}, {\bf 3}, 217-228.

\item Liptser, R.S. and Shiryayev, A.N. (1989) {\it The Theory of Martingales}, Kluwer, Dordrecht.

\item Mendy, I. (2013) Parametric estimation for sub-fractional Ornstein-Uhlenbeck process, {\it J. Stat. Plan. Infer.}, {\bf 143}, 663-674.

\item Michael, R. and Pfanzagl, J. (1971) The accuracy of the normal approximation for minimum contrast estimate, {\it Z.
Wahr. verw Gebeite}, 18:73-84.

\item Prakasa Rao, B.L.S. (1987) {\it Asymptotic Theory of Statistical Inference}, Wiley, New York.

\item Prakasa Rao, B.L.S. (1999) {\it Semimartingales and Their Statistical Inference}, CRC Press, Boca Raton and Chapman and Hall, London.

\item Prakasa Rao, B.L.S. (2010) {\it Statistical Inference for Fractional Diffusion Processes}, Wiley, London.

\item Prakasa Rao, B.L.S. (2017a) On some maximal and integral inequalities for sub-fractional Brownian motion, {\it Stochastic Anal. Appl.}, {\bf 35}, 279-287.

\item Prakasa Rao, B.L.S. (2017b) Optimal estimation of a signal perturbed by a sub-fractional Brownian motion, {\it Stochastic Anal. Appl.}, {\bf 35}, 533-541.

\item Prakasa Rao, B.L.S. (2017c) Parameter estimation for linear stochastic differential equations driven by sub-fractional Brownian motion, {\it Random Oper. and Stoch. Equ.}, {\bf 25}, 235-247.

\item Shen, G.J. and Yan, L.T. (2014) Estimators for the drift of subfractional Brownian motion, {\it Communications in Statistics-Theory and Methods}, {\bf 43}, 1601-1612.

\item Tudor, Constantin (2007a) Some properties of the sub-fractional Brownian motion, {\it Stochastics}, {\bf 79}, 431-448.

\item Tudor, Constantin (2007b) Prediction and linear filtering with sub-fractional Brownian motion, Preprint.

\item Tudor, Constantin (2008) Some aspects of stochastic calculus  for the sub-fractional Brownian motion, {\it Analele Universitat ii Bucaresti, Matematica}, Anul LVII, pp. 199-230.

\item Tudor, Constantin (2009) On the Wiener integral with respect to a sub-fractional Brownian motion on an interval, {\it J. Math. Anal. Appl.}, {\bf 351}, 456-468.

\item Yan, L., Shen, G. and He, K. (2011) Ito's formula for a sub-fractional Brownian motion, {\it Communications of Stochastic Analysis}, {\bf 5}, 135-159.

\end{description}
\end{document}